\documentclass[a4paper,12pt]{scrartcl} 
\usepackage[T1]{fontenc}
\usepackage{latexsym}
\usepackage{amssymb}
\usepackage{graphics}
\newtheorem{thm}{Theorem}

\newcommand{\C}{\mathbb{C}}

\newcommand{\E}{\mathbb{E}}
\newcommand{\lb}{\lambda}

\title{\textbf{On the Relative Weak Asymptotic Homomorphism Property for Triples of Group von Neumann Algebras}}
\author{Paul Jolissaint}

\begin{document}

\maketitle

\begin{abstract}
A triple of finite von Neumann algebras $B\subset N\subset M$ is said to have the relative weak asymptotic homomorphism property if there exists a net of unitaries $(u_i)_{i\in I}\subset U(B)$ such that
$$
\lim_{i\in I}\Vert \E_B(xu_iy)-\E_B(\E_N(x)u_i\E_N(y))\Vert_2=0
$$
for all $x,y\in M$.
Then recently, J. Fang, M. Gao and R. Smith proved that the triple $B\subset N\subset M$ has the relative weak asymptotic homomorphism property if and only if $N$ contains the set of all $x\in M$ such that $Bx\subset\sum_{i=1}^n x_iB$ for finitely many elements $x_1,\ldots,x_n\in M$. Furthermore, if 
$H<G$ is a pair of groups, they get a purely algebraic characterization of the weak asymptotic homomorphism property for the pair of von Neumann algebras $L(H)\subset L(G)$, but their proof requires a result which is very general and whose proof is rather long. We extend the result to the case of a triple of groups $H<K<G$, we present a direct and elementary proof of the above mentioned characterization and we introduce three more equivalent combinatorial conditions on the triple $H<K<G$, one of them stating that the subspace of $H$-compact vectors 
of the quasi-regular representation of $H$ on $\ell^2(G/H)$ is contained in $\ell^2(K/H)$.
\end{abstract}

\noindent
\emph{Key words:} von Neumann algebra, one sided quasi-normalizer, discrete group, quasi-regular representation, asymptotic homomorphism
\newline
\emph{AMS classification:} 46L10, 22D25

\section{Introduction}

Let $1\in B\subset N\subset M$ be a triple of finite von Neumann algebras gifted with a fixed, normal, finite, faithful and normalized trace $\tau$. Then $\E_N$ (resp. $\E_B$) denotes the $\tau$-preserving conditional expectation from $M$ onto $N$ (resp. $B$); we also set
$M\ominus N=\{x\in M: \E_N(x)=0\}$.
\par\vspace{3mm}

Following \cite{FGS}, we say that the triple $B\subset N\subset M$ has the \emph{relative weak asymptotic homomorphism property} if there exists a net of unitaries $(u_i)_{i\in I}\subset U(B)$ such that, for all $x,y\in M\ominus N$,
$$
\lim_{i\in I}\Vert \E_B(xu_i y)\Vert_2=0.
$$

The \emph{one sided quasi-normalizer of $B$ in $M$} is the set of elements $x\in M$ for which there exist finitely many elements $x_1,\ldots,x_n\in M$ such that $Bx\subset \sum_{i=1}^n x_iB$. It is denoted by $q\mathcal{N}_M^{(1)}(B)$.
\par\vspace{3mm}
Inspired by \cite{Chi}, the authors of \cite{FGS} prove in Theorem 3.1 that the triple $B\subset N\subset M$ has the relative weak asymptotic homomorphism property if and only if $q\mathcal{N}_M^{(1)}(B)\subset N$. Furthermore, they also study the case of group algebras that we recall now.

\par\vspace{3mm}
Let $G$ be a discrete group and let $H$ be a subgroup of $G$. Then there is a natural analogue of the one sided quasi-normalizer for such a pair of groups: we denote by $q\mathcal{N}_G^{(1)}(H)$ the set of elements $g\in G$ for which there exist finitely many elements $g_1,\ldots,g_n\in G$ such that $Hg\subset\cup_{i=1}^n g_iH$. 

\par\vspace{3mm}
Thus, if $H<K<G$ is a triple of groups, if $B=L(H)\subset N=L(K)\subset M=L(G)$ denotes the triple of von Neumann algebras associated to $H<K<G$, it is reasonable to ask whether $B\subset N\subset M$ has the relative weak asymptotic homomorphism property if and only if $q\mathcal{N}_G^{(1)}(H)\subset K$. Corollary 5.4 in \cite{FGS} states that this is indeed true when $K=H$, but the proof presented there relies heavily on the main theorem of the article. It is thus natural to look for a more direct and elementary proof of the above mentionned result, and the aim of the present note is to provide such a proof and to add three more equivalent conditions.


\section{The main result}
Before stating our result, let us fix some additional notations.
For each element $g\in G$ we denote by $\lambda_g$ the unitary operator acting by left translation on $\ell^2(G)$, i.e. $(\lambda_g\xi)(g')=\xi(g^{-1}g')$ for every $\xi\in\ell^2(G)$ and every $g'\in G$. We denote also by $L_f(G)$ the subalgebra of all elements of $L(G)$ with finite support, i.e. $L_f(G)$ is the linear span of $\lambda(G)$ in $B(\ell^2(G))$. 
\par\vspace{3mm}
We fix a triple of groups $H<K<G$ for the rest of the article.
\par\vspace{3mm}
Let $\pi$ denote the quasi-regular representation of $G$ on $\ell^2(G/H)$; we denote by $[g]$ the equivalence class $[g]=gH$, so that $\pi(g)\xi([g'])=\xi([g^{-1}g'])$ for all $g,g'\in G$ and $\xi\in\ell^2(G/H)$. Following \cite{BR}, we say that a vector $\xi\in\ell^2(G/H)$ is $H$-\emph{compact} if the norm closure of its $H$-orbit $\{\pi(h)\xi:h\in H\}$ is a compact subset of $\ell^2(G/H)$. The set of all $H$-compact vectors is a closed subspace of $\ell^2(G/H)$ that we denote by $\ell^2(G/H)_{c,H}$.
We also set
$$
\ell^2(G/H)^H=\{\xi\in\ell^2(G/H): \pi(h)\xi=\xi\ \forall h\in H\},
$$
which is the subspace of all $H$-invariant vectors of $\ell^2(G/H)$. We observe that it is contained in $\ell^2(G/H)_{c,H}$.

\begin{thm}
Let $H<K<G$ and $B=L(H)\subset N=L(K)\subset M=L(G)$ be as above. Then the following conditions are equivalent:
\begin{enumerate}
	\item [(1)] There exists a net $(h_i)_{i\in I}\subset H$ such that, for all $x,y\in M\ominus N$, one has
	$$
	\lim_{i\in I}\Vert \E_B(x\lb_{h_i} y)\Vert_2=0,
$$
i.e. the net of unitaries in the relative weak asymptotic homomorphism property may be chosen in the subgroup $\lb(H)$ of $U(B)$.
	\item [(2)] The triple $B\subset N\subset M$ has the relative weak asymptotic homomorphism property.
	\item [(3)] If $g\in G$ and $F\subset G$ finite are such that $Hg\subset FH$, then $g\in K$, i.e. $q\mathcal{N}_G^{(1)}(H)\subset K$.
	\item [(4)] The subspace of $H$-compact vectors $\ell^2(G/H)_{c,H}$ is contained in $\ell^2(K/H)$.
	\item [(5)] The subspace $\ell^2(G/H)^H$ is contained in $\ell^2(K/H)$.
	\item [(6)] For every non empty finite set $F\subset G\setminus H$, there exists $h\in H$ such that 
	$$
	FhF\cap H=\emptyset.
	$$
\end{enumerate}
\end{thm}
\emph{Proof.} \emph{(1) $\Rightarrow$ (2)} is obvious.
\par\vspace{1mm}\noindent
\emph{(2) $\Rightarrow$ (3).}
 Observe that condition (3) is equivalent to the following statement (since, if $g\notin K$, then $HgH\cap K=\emptyset$): 
\par
\emph{For every $g\in G\setminus K$, and for every non empty finite set $F\subset G\setminus K$, there exists $h\in H$ such that $Fhg\cap H=\emptyset$.} 
\par
Thus, let us assume that condition (3) does not hold. There exists $g\in G\setminus K$ and a non empty finite set $F\subset G\setminus K$ such that $Fhg\cap H\not=\emptyset$ for every $h\in H$. Then let $u\in U(B)$. One has:
\begin{eqnarray*}
\sum_{g'\in F}\Vert \E_B(\lambda_{g'}u\lambda_g)\Vert_2^2 
&=&
\sum_{g'\in F}\left(\sum_{h\in H,g'hg\in H}|u(h)|^2\right)\\
&=&
\sum_{h\in H}\left(\sum_{g'\in F,g'hg\in H}|u(h)|^2\right)\\
&\geq &
\sum_{h\in H}|u(h)|^2=\Vert u\Vert_2^2=1
\end{eqnarray*}
since, for every $h\in H$, one can find $g'(h)\in F$ such that $g'(h)hg\in H$. Hence there cannot exist a net $(u_i)_{i\in I}\subset U(B)$ as above, and the triple $B\subset N\subset M$ does not have the relative weak asymptotic homomorphism property.
\par\vspace{1mm}\noindent
\emph{(3) $\Rightarrow$ (4).}
We choose a set of representatives $T\ni e$ of left classes so that $G=\sqcup_{t\in T}tH$, and let $\xi\not=0$ be an $H$-compact vector.
\par
Let $s\in T$ be such that $\epsilon:=|\xi([s])|>0$. There exists then finitely many vectors $\xi_1,\ldots,\xi_n\in\ell^2(G/H)$ such that, for every $h\in H$, there exists $1\leq j\leq n$ such that $\Vert\pi(h)\xi-\xi_j\Vert\leq\epsilon/2$. Set
$$
F=\bigcup_{j=1}^n\{t\in T:|\xi_j([t])|\geq\epsilon/2\},
$$
which is a finite set. Then we claim that $Hs\subset FH$. Indeed, if $h\in H$, let $t\in T$ be such that $[hs]=[t]$, and let $j$ be such that $\Vert\pi(h)\xi-\xi_j\Vert\leq\epsilon/2$. Then
$$
\epsilon-|\xi_j([t])|=|\xi([s])|-|\xi_j([t])|\leq |\xi([s])-\xi_j([hs])|\leq 
\Vert\pi(h)\xi-\xi_j\Vert\leq\epsilon/2
$$
hence $\epsilon/2\leq |\xi_j([t])|$ and $t\in F$. Thus $Hs\subset FH$, and condition (3) implies that $s\in K$. This proves that $\xi\in\ell^2(K/H)$.
\par\vspace{1mm}\noindent
\emph{(4) $\Rightarrow$ (5)} is obvious.
\par\vspace{1mm}\noindent
\emph{(5) $\Rightarrow$ (6).} Let us assume that the triple $H<K<G$ satisfies condition (5) but not (6). Then there exists a finite set $F=F^{-1}\subset G\setminus K$ such that $FhF\cap H\not=\emptyset$ for every $h\in H$. Set 
$$
\xi=\sum_{g\in F}\delta_{[g]}.
$$
Then $\xi\perp\ell^2(K/H)$, and one has for every $h\in H$:
$$
\langle\pi(h)\xi,\xi\rangle=\sum_{g,g'\in F}\langle\delta_{[hg]},\delta_{[g']}\rangle\geq 1
$$
since the condition on $F$ implies that for every $h\in H$, there exist $g,g'\in F$ such that $hgH=g'H$. Let $C$ be the closed convex hull of $\{\pi(h)\xi:h\in H\}$. Then it is easy to see that $\langle \zeta,\xi\rangle\geq 1$ for every $\zeta\in C$. Let $\eta\in C$ be the vector with minimal norm. By its uniqueness, it is $H$-invariant and non zero by the above observation. Thus, $\eta$ is supported in $K/H$ and orthogonal to $\ell^2(K/H)$ since $\xi$ is.
This is the expected contradiction.
\par\vspace{1mm}\noindent
\emph{(6) $\Rightarrow$ (1).} Let $I=\{F \subset G\setminus K: F\not=\emptyset,\ \textrm{finite}\}$ be the directed set of all non empty finite subsets of $G\setminus K$. Condition (6) states that, for every $F\in I$, there exists $h_F\in H$ such that $Fh_FF\cap H=\emptyset$. Let $x$ and $y$ in $L_f(G)$ that satisfy $\E_N(x)=\E_N(y)=0$. Let then $F_0\in I$ be so that the supports of $x$ and $y$ are contained in $F_{0}$. Then $x\lambda_{h_F}y=\sum_{g,g'\in F_{0}}x(g)y(g')\lambda_{gh_Fg'}$ for every $F\supset F_0$, thus $\E_B(x\lambda_{h_F}y)=0$ for every $F\supset F_0$. This proves that the triple $B\subset N\subset M$ satisfies condition (1) by density of $L_f(G)$ in $L(G)$.
\hfill
$\square$

\par\vspace{3mm}\noindent
\textbf{Remark.} In the case of a pair of groups $H<G$, which corresponds to $H=K$, condition (4) means that all $H$-invariant vectors in $\ell^2(G/H)$ are multiples of $\delta_{[e]}$, and this means that the unitary representation $\rho$ of $H$ on the subspace $\ell^2(G/H)\ominus\C\delta_{[e]}$ is ergodic in the sence of \cite{BR}.

\par
\bibliographystyle{plain}
\bibliography{max}

\vspace{1cm}
\noindent
\begin{flushright}
     \begin{tabular}{l}
       Universit\'e de Neuch\^atel,\\
       Institut de Math\'emathiques,\\       
       Emile-Argand 11\\
       CH-2000 Neuch\^atel, Switzerland\\
       \small {paul.jolissaint@unine.ch}
     \end{tabular}
\end{flushright}

\end{document}